\documentclass[11pt]{article}
\usepackage{latexsym}
\usepackage{amsmath}
\usepackage{amsfonts}
\usepackage{mathrsfs}

\newtheorem{lem}{Lemma}[section]
\newtheorem{thm}{Theorem}[section]

\newtheorem{defn}{Definition}[section]

\newcommand{\Rs}{\mathbb{R}}

\newcommand{\mm}{\bar{m} }

\newcommand{\rr}{\bar{r} }

\newcommand{\EE}{{\cal E} }

\newcommand{\bpr}{{\bf Proof.} \hspace{1 em}}
\newcommand{\epr}{ \\ \hspace*{4.5in} $\Box$ }
\newcommand{\beq}{ \begin{equation} }
\newcommand{\eeq}{ \end{equation} }
\newcommand{\bt}{ \begin{tabular} }
\newcommand{\et}{ \end{tabular} }

\begin{document}

\bibliographystyle{plain}
\title{On Affine Motions and Bar Frameworks in General Position} 
\vspace{0.3in}           
        \author{ A. Y. Alfakih 
  \thanks{E-mail: alfakih@uwindsor.ca} 
  \thanks{Research supported by the Natural Sciences and Engineering
         Research Council of Canada.} 
\\
          Department of Mathematics and Statistics \\
          University of Windsor \\
          Windsor, Ontario N9B 3P4 \\
          Canada  
\and
               Yinyu Ye
\thanks{E-mail:yinyu-ye@stanford.edu}
\thanks{Research supported in part by NSF Grant GOALI 0800151 
                   and DOE Grant DE-SC0002009.    }
\\
         Department of Management Science and Engineering \\
         Stanford University \\
         Stanford, California 94305 \\
         USA }

\date{\today}  
\maketitle

\noindent {\bf AMS classification:}  52C25, 05C62.  

\noindent {\bf Keywords:} Bar frameworks,   
universal rigidity, stress matrices, points in general position, Gale transform.  
\vspace{0.1in}

\begin{abstract}
A configuration $p$ in $r$-dimensional Euclidean space is a finite  
collection of points $(p^1,\ldots, p^n)$ that affinely span $\Rs^r$.  
A bar framework, denoted by $G(p)$, in $\Rs^r$ is a simple 
graph $G$ on $n$ vertices together with 
a configuration $p$ in $\Rs^r$. A given bar framework $G(p)$ is said 
to be universally rigid if
there does not exist another configuration $q$ in any Euclidean space, 
not obtained from $p$ by a rigid motion, 
such that $||q^i-q^j||$ = $||p^i-p^j||$  for each edge $(i,j)$ of $G$. 

It is known \cite{alf10, con99} that if configuration $p$ is generic and
bar framework $G(p)$ in $\Rs^r$ admits a positive semidefinite
stress matrix $S$ of rank ( $n-r-1$), then $G(p)$ is universally rigid. 
Connelly asked \cite{con09} whether the same result holds true if the
genericity assumption of $p$ is replaced by the weaker assumption of
general position. We answer this question in the affirmative in this paper. 
\end{abstract}

\section{Introduction}

A {\em configuration} $p$ in $r$-dimensional Euclidean space is a   
finite collection of points $(p^1,\ldots,p^n)$ in $\Rs^r$ 
that affinely span $\Rs^r$.  
A {\em bar framework} (or framework for short) in $\Rs^r$, 
denoted by $G(p)$, is  
a configuration $p$ in $\Rs^r$ together with a simple graph $G$
on the vertices $1,2, \ldots,n$. For a simple graph $G$, we denote its 
node set by $V(G)$ and its edge set by $E(G)$. To avoid trivialities,  
we assume throughout this paper that graph $G$ is connected and not complete.

Framework $G(q)$ in $\Rs^r$ is said to be  
{\em congruent} to framework $G(p)$ in $\Rs^r$ if configuration $q$ is 
obtained from configuration $p$ by a rigid motion. That is,   
if $|\!|q^i-q^j|\!|$= $|\!|p^i-p^j|\!|$ for all
$i,j=1,\ldots,n$, where $|\!|.|\!|$ denotes the Euclidean norm. 
We say that framework $G(q)$ in $\Rs^s$ is {\em equivalent } to framework
$G(p)$ in $\Rs^r$ if $|\!|q^i-q^j|\!|$= $|\!|p^i-p^j|\!|$ for all
$(i,j) \in E(G)$.  
Furthermore, 
we say that framework $G(q)$ in $\Rs^r$ is {\em affinely-equivalent } 
to framework $G(p)$ in $\Rs^r$ if $G(q)$ is equivalent to $G(p)$ and
configuration $q$ is obtained from configuration $p$ by an affine motion;
i.e., $q^i = A p^i + b$, for all $i=1,\ldots,n$, 
for some $r \times r$ matrix $A$ and an $r$-vector $b$.  

A framework $G(p)$ in $\Rs^r$ is said to be {\em universally rigid}
if there does exist a framework $G(q)$ in any Euclidean space that
is equivalent, but not congruent, to $G(p)$. 
The notion of a stress matrix $S$ of a framework $G(p)$ plays a 
key role in the problem of universal rigidity of $G(p)$.

\subsection{Stress Matrices and Universal Rigidity} 

Let $G(p)$ be a framework on $n$ vertices in $\Rs^r$.
An {\em equilibrium stress} of $G(p)$ is a real valued function
$\omega$ on $E(G)$ such that
\beq \label{seq}  
\sum_{j:(i,j) \in E(G)} \omega_{ij} (p^i - p^j) =0 \mbox{ for all } i=1,\ldots,n.
\eeq 
 
Let $\omega$ be an equilibrium stress of $G(p)$. Then the 
$n \times n$ symmetric matrix $S=(s_{ij})$ where
\beq \label{defS} 
s_{ij} = \left\{ \begin{array}{ll} -\omega_{ij} & \mbox{if } (i,j) \in E(G), \\
                        0   & \mbox{if } i \neq j \mbox{ and } (i,j) \not \in E(G), \\
                   {\displaystyle \sum_{k:(i,k) \in E(G)} \omega_{ik}} & \mbox{if } i=j, 
                   \end{array} \right. 
\eeq
is called the {\em stress matrix} associated with $\omega$, or a stress matrix
of $G(p)$.  
The following result provides a sufficient condition for the 
universal rigidity of a given framework. 

\begin{thm}[Connelly \cite{con82,con99}, Alfakih \cite{alf07a}] \label{thma} 
Let $G(p)$ be a bar framework in $\Rs^r$, for some $r \leq n-2$. If the 
following two conditions hold: 
\begin{enumerate}
\item There exists a positive semidefinite stress matrix $S$ of $G(p)$ of rank
($n-r-1$).  
\item There does not exist a bar framework $G(q)$ in $\Rs^r$
that is affinely-equivalent, but not congruent, to $G(p)$. 
\end{enumerate}
Then $G(p)$ is universally rigid. 
\end{thm}

Note that ($n-r-1$) is the maximum possible value for the rank of the
stress matrix $S$.  
In connection with Theorem \ref{thma}, we mention the following result
obtained in So and Ye \cite{SY05} and Biswas et al. \cite{bty08}: 
Given a framework $G(p)$ in $\Rs^r$, 
if there does not exist a framework $G(q)$ in $\Rs^s$ ( $s \neq r$) that 
is equivalent to $G(p)$, then $G(p)$ is universally rigid. Moreover, if $G(p)$ 
contains a clique of $r+1$ points in general position, then the existence of
a rank-($n-r-1$) positive semidefinite stress matrix implies that 
framework $G(p)$ is universally rigid, regardless whether the other 
non-clique points are in general position or not. 

Condition 2 of Theorem \ref{thma} is satisfied if configuration
$p$ is assumed to be generic (see Lemma \ref{conicpge} below). 
A configuration $p$ (or a framework $G(p)$) is said to be {\em generic} if 
all the coordinates of $p^1,\ldots,p^n$ are algebraically independent 
over the integers. That is, if there does not exist
a non-zero polynomial $f$ with integer coefficients 
such that $f(p^1,\ldots,p^n)=0$. Thus  

\begin{thm}[Connelly \cite{con99}, Alfakih \cite{alf10}] \label{thmac} 
Let $G(p)$ be a generic bar framework on $n$ nodes in $\Rs^r$, for some $r \leq n-2$. 
If there exists a positive semidefinite stress matrix $S$ of $G(p)$ of rank
($n-r-1$). Then $G(p)$ is universally rigid. 
\end{thm}

The converse of Theorem \ref{thmac} is also true.

\begin{thm}[Gortler and Thurston \cite{gt09}] 
Let $G(p)$ be a generic bar framework on $n$ nodes in $\Rs^r$, for some $r \leq n-2$. 
If $G(p)$ is universally rigid, then  
there exists a positive semidefinite stress matrix $S$ of $G(p)$ of rank
($n-r-1$).  
\end{thm}

Connelly \cite{con09} asked whether a result similar to 
Theorem \ref{thmac} holds if the genericity
assumption of $G(p)$ is replaced by the weaker assumption of general position.  
A configuration $p$ (or a framework $G(p)$) in $\Rs^r$ is 
said to be in {\em general position} if 
no subset of the points $p^1,\ldots,p^n$ 
of cardinality $r+1$ is affinely dependent.
For example, a set of points in the plane are in general position
if no 3 of them lie on a straight line. 

In this paper we answer Connelly's question in the affirmative. Thus the
following theorem is the main result of this paper. 

\begin{thm} \label{main} 
Let $G(p)$ be a bar framework on $n$ nodes in general position in $\Rs^r$, 
for some $r \leq n-2$. 
If there exists a positive semidefinite stress matrix $S$ of $G(p)$ of rank
($n-r-1$). Then $G(p)$ is universally rigid. 
\end{thm}

The proof of Theorem \ref{main} will be given in Section \ref{pr}. 
This paper and \cite{aty10c} are first steps toward the study of universal
rigidity under the general position assumption. In \cite{aty10c}, it was
shown that the framework $G(p)$ on $n$ nodes in general position in $\Rs^r$ for
some $r \leq n-2$, where $G$ is the $(r+1)$-lateration graph, 
admits a rank $(n-r-1)$ positive semi-definite stress matrix. 
 
\section{Preliminaries} 

To develop the ingredients needed for the proof of our main result, 
we review the necessary background on affine motions, stress matrices,
and Gale matrices. 
 
An affine motion in $\Rs^r$ is a map $f: \Rs^r \rightarrow \Rs^r$ of
the form 
\[
f(p^i)= A p^i + b,
\]  
for all $p^i$ in $\Rs^r$, where $A$ is 
an $r \times r$ matrix and $b$ is an $r$-vector. 
A rigid motion is an affine motion where matrix $A$ is orthogonal.  

Vectors $v^1,\ldots,v^m$ in $\Rs^r$ are said to lie on a {\em quadratic
at infinity} if there exists a non-zero symmetric $r \times r$ 
matrix $\Phi$ such that 
\beq
({v^i})^T \Phi v^i = 0 , \mbox{ for all } i=1,\ldots,m.
\eeq
 
\begin{lem}(Connelly \cite{con05}) \label{conicp} 
Let $G(p)$ be a bar framework on $n$ vertices in $\Rs^r$. Then the following
two conditions are equivalent:
\begin{enumerate}
\item There exists
a framework $G(q)$ in $\Rs^r$ that is equivalent, but not congruent, 
to $G(p)$ such that ${q^i} = A p^i +b$ for all $i=1,\ldots,n$, 
\item The vectors $p^i-p^j$ for all $(i,j) \in E(G)$ lie on a quadratic at infinity.
\end{enumerate}
\end{lem}

\begin{lem}(Connelly \cite{con05}) \label{conicpge} 
Let $G(p)$ be a generic bar framework on $n$ vertices in $\Rs^r$. Assume that
each node of $G$ has degree at least $r$. Then  
the vectors $p^i-p^j$ for all $(i,j) \in E(G)$ do not lie on a quadratic at infinity.
\end{lem}

Therefore, under the genericity assumption,  
Condition 2 in Lemma \ref{conicp} does not hold. Consequently, Theorem 
\ref{thmac} follows as a simple corollary of Theorem \ref{thma}.  
  
Note that Condition 2 in Lemma \ref{conicp} is expressed in terms of the
edges of $G$. An equivalent condition in terms of the missing edges of 
$G$ can also be obtained using Gale matrices. This equivalent condition 
turns out to be crucial for our proof of Theorem \ref{main}. 

To this end,  
let $G(p)$ be a framework on $n$ vertices in $\Rs^r$. Then  
the following $(r+1) \times n$ matrix 
\beq \label{defA} 
{\cal A} := \left[ \begin{array}{cccc} p^1 & p^2 & \ldots & p^n \\ 
                           1   & 1 & \ldots & 1 
         \end{array}  \right] 
\eeq 
has full row rank since $p^1,\ldots,p^n$ affinely span $\Rs^r$. 
Note that $r \leq n-1$. Let 
\beq \label{defrr} 
\rr = \mbox{the dimension of the null space
of } {\cal A}; \mbox{ i.e., } \rr=n-1-r.  
\eeq
\begin{defn} \label{defLam} 
Suppose that the null space of ${\cal A}$ is nontrivial, i.e., $\rr \geq 1$.
Any $n \times \rr$ matrix $Z$ 
whose columns form a basis of the null space of ${\cal A}$    
is called a {\em Gale matrix} of configuration $p$.  
Furthermore, the $i$th row of $Z$, considered as
a vector in $\Rs^{\rr}$, is called a {\em Gale transform} of $p^i$ \cite{gal56}. 
\end{defn} 

Let $S$ be a stress matrix of $G(p)$ then 
it follows from (\ref{defS}) and (\ref{defA}) that 
\beq \label{SZ1} 
{\cal A} S =0. 
\eeq
Thus
\beq \label{SZ2} 
S = Z \Psi Z^T ,
\eeq
for some $\rr \times \rr$ symmetric matrix $\Psi$, where $Z$ is a
Gale matrix of $p$. It immediately follows from (\ref{SZ2}) that 
rank $S$ = rank $\Psi$. Thus, $S$ attains its maximum rank of $\rr=(n-1-r)$ if and
only if $\Psi$ is nonsingular, i.e., rank $\Psi = \rr$. 
 
Let $e$ denote the vector of all 1's in $\Rs^n$, and  
let $V$ be an $n \times (n-1)$ matrix that satisfies: 
\beq \label{defV}
V^Te=0, \; V^TV=I_{n-1}, 
\eeq
where $I_{n-1}$ is the identity matrix of order $(n-1)$. 
Further, let $E^{ij}$, $i \neq j$, denote the $n \times n$ symmetric matrix with 1 in
the $(i,j)$th and $(j,i)$th entries and zeros elsewhere, and let
$\EE(y) = \sum_{(i,j) \not \in E(G)} y_{ij} E^{ij} $ where $y_{ij}=y_{ji}$. 
In other words, the $(k,l)$ entry of matrix $\EE(y)$ is given by  
\beq \label{defEE} 
 \EE(y)_{kl} = \left\{ \begin{array}{ll} 0 & \mbox{ if } (k,l) \in E(G), \\
                                   0 & \mbox{ if } k=l, \\
                         y_{kl} & \mbox{ if } k \neq l \mbox{ and } (k,l) \not \in E(G). 
    \end{array} \right.
\eeq
 
Then we have the following result.

\begin{lem}(Alfakih \cite{alf10}) \label{lll} 
Let $G(p)$ be a bar framework on $n$ vertices in $\Rs^r$ and let $Z$ be any 
Gale matrix of $p$. Then 
the following two conditions are equivalent:
\begin{enumerate}
\item The vectors $p^i-p^j$ for all $(i,j) \in E(G)$ lie on a quadratic at infinity.
\item There exists a non-zero $y=(y_{ij}) \in \Rs^{\mm}$ such that: 
\beq \label{conicZ} 
V^T \EE(y) Z = {\bf 0}, 
\eeq 
where $\mm$ is the number of missing edges of $G$, $V$ is defined in (\ref{defV}),
and $\EE(y)$ is defined in (\ref{defEE}). {\bf 0} here is the zero matrix of
dimension $(n-1) \times \rr$. 
\end{enumerate} 
\end{lem}

Condition 2 of Lemma \ref{lll} can be easily understood if 
a projected Gram matrix approach is used for   
the universal rigidity of bar frameworks (see \cite{alf10} for details).  

\section{Proof of Theorem \ref{main}} \label{pr} 

The main idea of the proof is to show that Condition 2 of Lemma \ref{lll}
does not
hold under the general position assumption, and under the assumption that
$G(p)$ admits a positive semidefinite stress matrix of rank ($n-r-1$).    
The choice of the particular Gale matrix to be used in equation (\ref{conicZ}) 
is critical in this regard.  We begin with a few necessary lemmas.
 
\begin{lem} \label{Zli} 
Let $G(p)$ be a framework on $n$ nodes in general position in $\Rs^r$ and
let $Z$ be any Gale matrix of configuration $p$.
Then any $\rr \times \rr$ submatrix  of $Z$ is nonsingular. 
\end{lem}

\bpr For a proof see e.g., \cite{alf07a}. 
\epr

Let $\bar{N}(i)$ denote the set of nodes of graph $G$ that are non-adjacent 
to node $i$; i.e., 
\beq
\label{defNb}  \bar{N}(i) = \{ j \in V(G): j\neq i \mbox{ and } (i,j) \not \in E(G) \}, 
\eeq 

\begin{lem} \label{lemZh} 
Let $G(p)$ be a framework on $n$ nodes in general position in $\Rs^r$. 
Assume that $G(p)$ has a stress matrix $S$ of rank
$( n-1-r) $. Then there exists a Gale matrix $\hat{Z}$ of $G(p)$ such that  
$\hat{z}_{ij}=0$ for all $j=1,\ldots,\rr$ and $i \in \bar{N}(j+r+1)$. 
\end{lem}

\bpr 
Let $G(p)$ be in general position in $\Rs^r$ and assume that it has
a stress matrix $S$ of rank $\rr= (n-1-r )$. 
Let $Z$ be any Gale matrix of $G(p)$, then 
$S = Z \Psi Z^T$ for some non-singular symmetric 
$\rr \times \rr$ matrix $\Psi$.   
Let us write $Z$ as: 
\beq
Z = \left[ \begin{array}{c} Z_1 \\ Z_2 \end{array} \right],
\eeq
where $Z_2$ is $\rr \times \rr$. By Lemma  \ref{Zli}, $Z_2$
is non-singular. Now let   
\beq \label{defZh} 
\hat{Z}=(\hat{z}_{ij})= Z \Psi {Z_2}^T. 
\eeq 
Then $\hat{Z}$ is a Gale matrix of $G(p)$. This  
simply follows from the fact that the matrix obtained by multiplying  
any Gale matrix of $G(p)$ from the right by a non-singular $\rr \times \rr$ 
matrix, is also a Gale matrix of $G(p)$. Furthermore,  

\[
S = Z \Psi Z^T = Z \Psi \, [ Z_1^T \;\; Z_2^T] = [Z \Psi Z_1^T \;\; \hat{Z}].  
\]
In other words, $\hat{Z}$ consists of the last $\rr$ columns of $S$.
Thus $\hat{z}_{ij} = s_{i,j+r+1}$. By the definition of $S$ we have 
$s_{i,j+r+1}=0$ for all $i \neq j+r+1$ and $(i,j+r+1) \not \in E(G)$. Therefore, 
$\hat{z}_{ij} = 0$ for all $j=1,\ldots,\rr$ and $i \in \bar{N}(j+r+1)$. 
\epr

\begin{lem} \label{lemsys2} 
Let the Gale matrix in (\ref{conicZ}) be $\hat{Z}$ as defined in (\ref{defZh}). 
Then the system of equations (\ref{conicZ}) is equivalent to  the system  
of equations
\beq \label{conicZ2}
\EE(y) \hat{Z} = {\bf 0}. 
\eeq
{\bf 0} here is the zero matrix of dimension $n \times \rr$. 
\end{lem}

\bpr
System of equations (\ref{conicZ}) is equivalent to the following system of 
equations in the unknowns, $y_{ij}$ ($i \neq j$ and $(i,j) \not \in E(G)$) 
and $\xi \in \Rs^{\rr}$: 
\beq \label{eq1} 
 \EE(y) \hat{Z} = e \, \xi^T,
\eeq 
Now for $j=1,\ldots,\rr$, we have that the $(j+r+1,j)$th entry of
$\EE(y) \hat{Z}$ is equal to $\xi_j$. But using (\ref{defEE}) and Lemma \ref{lemZh}
we have  
\[ 
(\EE(y) \hat{Z})_{j+r+1 ,j}  =  \sum_{i=1}^n \EE(y)_{j+r+1, i} \; \hat{z}_{ij}   
                   = \sum_{i: i \in \bar{N}(j+r+1)} y_{j+r+1, i} \; \hat{z}_{ij}  =0. 
\]
Thus, $\xi = 0$ and the result follows.
\epr  

Now we are ready to prove our main theorem.

\noindent{\bf Proof of Theorem \ref{main} } 

Let $G(p)$ be a framework on $n$ nodes in general position in $\Rs^r$. 
Assume that $G(p)$ has a positive semidefinite stress matrix $S$ of rank
$\rr= n-1-r $. Then 
deg$(i) \geq r+1$ for all $i \in V(G)$, 
i.e., every node of $G$ is adjacent to at least $r+1$ nodes 
(for a proof see \cite[Theorem 3.2]{alf07a}). Thus   
\beq \label{deg}  
|\bar{N}(i)| \leq n-r-2 = \rr - 1.  
\eeq

Furthermore, 
it follows from  Lemmas \ref{lemZh}, \ref{lemsys2}  
and \ref{lll} that 
the vectors $p^i-p^j$ for all $(i,j) \in E(G)$ lie on a quadratic at infinity
if and only if system of equations (\ref{conicZ2}) has a non-zero solution $y$.  
But (\ref{conicZ2}) can be written as  

\[
\sum_{j: \in \bar{N}(i)} y_{ij} \hat{z}^j = 0 , \mbox{ for } i=1,\ldots,n,
\]
where $(\hat{z}^i)^T$ is the $i$th row of $\hat{Z}$. Now  it follows from (\ref{deg})
that  $y_{ij}=0$ for all $(i,j) \not \in E(G)$ since by Lemma \ref{Zli}
any subset of $\{\hat{z}^1,\ldots,\hat{z}^n\}$ of cardinality $\leq \rr-1$ is linearly 
independent.

Thus system (\ref{conicZ2}) does not have a nonzero solution $y$. Hence
the vectors $p^i-p^j$, for all $(i,j) \in E(G)$, do not lie on a quadratic at infinity.
Therefore, by Lemma \ref{conicp}, there does not exist a framework $G(q)$ in $\Rs^r$
that is affinely-equivalent, but not congruent, to $G(p)$. Thus by Theorem \ref{thma}, 
$G(p)$ is universally rigid. 
\epr 


\end{document}